\input AHTOH-E.STY
\def\tphi{{\~{\raise5pt\hbox{$\scriptstyle\phi$}}}}

\UDC{
512.542        
+ 512.543.72   
}

\MSC{
20D60,       
20F70
}

\title{%
What do
\\
Frobenius's, Solomon's, and Iwasaki's theorems
on divisibility in groups
\\
have in common?
}
\author{%
Elena K. Brusyanskaya$^\sharp$
\quad
Anton A. Klyachko$^\sharp$
\quad
Andrey V. Vasil'ev$^{\flat\natural}$
}
\address{
$^\sharp$
Faculty of Mechanics and Mathematics of Moscow State University,
Moscow 119991, Leninskie gory, MSU.\\
\strut
$^\flat$Sobolev Institute of Mathematics, Novosibirsk 630090,
prospekt akademika Koptyuga 4.
\\
\strut
$^\natural$Novosibirsk State University,
Novosibirsk 630090, ul. Pirogova, 1.
\\
ebrusianskaia@gmail.com
\quad
klyachko@mech.math.msu.su
\quad
vasand@math.nsc.ru
}

\grantsFirstSecond{\RFBR19-01-00591}%
\grantsThird{%
the program of fundamental scientific researches of the SB RAS
no.I.1.1., project no.0314-2016-0001}

\abstract{%
Our result contains as special cases the
Frobenius theorem (1895) on the number of
solutions to the equation $x^n=1$
in a group,
the Solomon theorem (1969) on the
number of solutions in a group to a
system of equations having
fewer
equations than unknowns,
and
the Iwasaki theorem (1985)
on roots of subgroups.
There are other curious corollaries on
groups and rings.
}

\s 0.
Introduction

The following result was proved in XIX century.

\proclaim
Frobenius theorem {\rm[Frob95] (see also [And16])}.
The number of solutions to the equation $x^n=1$ in
a finite group
is divisible
by $\GCD(|G|,n)$ for any integer $n$.

This theorem was generalised in different directions,
see, e.g.,
[Hall36],
[Kula38],
[Sehg62],
[BrTh88],
[Yosh93],
[AsTa01],
[ACNT13],
and references therein.
For example, Frobenius himself [Frob03] obtained
the following generalisation in 1903:
\disp{\sl
for any positive integer $n$ and any element $g$ of a finite group
$G$, the number of solutions
to the equation $x^n=g$ in $G$ is divisible by
the
greatest common divisor of $n$ and the order of
the centraliser of $g$;
}
Ph. Hall
{\rm([Hall36], Theorem II)} showed that
\disp{\sl
\narrower
\narrower
\narrower
in any finite group, the number of solutions
to a system of equations in one
unknown is divisible
by $\GCD(|C|,n_1,n_2,\dots)$, where $C$ is
the
centraliser of the set of all coefficients and~$n_i$ are exponent sums
of the unknown in the $i$-th equation.
}%
Here, as usual, an equation over a group $G$
is an
expression of the form $v(x_1,\dots,x_m)=1$, where~$v$ is a word whose
letters are unknowns, their inverses, and
elements of~$G$ (called \emph{coefficients}).
In other terms, the left-hand side of an equation is an element
of the free
product $G*F(x_1,\dots,x_m)$ of $G$ and the free group
$F(x_1,\dots,x_m)$
of rank $m$ (where~$m$ is the number of unknowns).


The following theorem is also
about equations in groups and divisibility, but on the first
view, it is not similar to the Frobenius theorem and its
generalisations.

\proclaim
Solomon theorem \rm [Solo69].
In any group, the
number of solutions to a system of coefficient-free equations is
divisible by the order of the group provided the
number of
equations is less than the number of
unknowns.

This theorem was also generalised in different directions, see
[Isaa70],
[Stru95],
[AmV11],
[GRV12],
[KM14],
[KM17],
and references therein.
For instance, in [KM14], it was shown that
\disp{\sl
in any group,
the number of solutions to a system of equations (with coefficients
from this group)
is divisible by
the order of the intersection of centralisers of all coefficients
provided
the rank of the matrix composed of the exponent sums of 
the $j$-th
unknown in 
the $i$-th 
equation is less than the number of unknowns.
}
Solomon himself wrote in [Solo69]:
\disp{\it
``There seems to
be no connection between this theorem and the Frobenius theorem on
solutions of $x^k=1$.''
}
Nevertheless, a connection between the Frobenius and Solomon theorems
exists.

\Th 1%
\fn{%
{\bf Theorem 0} in the journal version.}%
.
In any (not necessarily finite) group, the number of solutions to a (not
necessarily finite) system of equations in $m$ unknowns is a multiple of
the greatest common divisor of the centraliser of the set of coefficients
and the number~$\Delta_m\over\Delta_{m-1}$, where $\Delta_i$ is the
greatest common divisor of all minors of order $i$ of the matrix of the
system, and the following conventions are assumed: $\Delta_i=0$ if $i$
larger than the number of equations, $\Delta_0=1$, and ${0\over0}=0$.

We define the
\emph{greatest common divisor $\GCD(G,n)$}
of a group $G$ and an integer $n$ as
the least common multiple of orders of subgroups of $G$ dividing
$n$.
The divisibility is always understood in the sense of cardinal arithmetic:
each infinite cardinal is divisible by all smaller nonzero cardinals
(and surely
zero is divisible by all cardinals and divides only
zero). This means that $\GCD(G,0)=|G|$ for any group $G$ and, e.g.,
$\GCD(\SL_2(\Z),2018)=2$.
Although, the reader will not lose much by
assuming all group to be finite; in this case,
$\GCD(G,n)=\GCD(|G|,n)$ by the Sylow theorem (and because a finite
$p$-group contains subgroups of all possible orders).

The
\emph{matrix of a system of equations over a group}
is the integer matrix
$A=(a_{ij})$, where $a_{ij}$ is
the exponent sum of 
the $j$-th 
unknown in
the $i$-th 
equation.
For example, the matrix of the system
$$
\cases{
xay^2[x,y]^{\the\year}(xby)^3=1\cr
bx^3y[x,y]^{100}(xby)^4=1\cr
[x,y^5]x^{-2}=1\cr
}
$$
(where $x$ and $y$ are unknowns, and $a$ and $b$ are coefficients, i.e.
some fixed group elements) has the form
$$
\pmatrix{
4&5\cr
7&5\cr
-2&0\cr
}.
$$
As usual,
the \emph{minors of order $i$} are determinants
of
submatrices composed of entries at the intersections of some $i$
rows and $i$ columns. In the example above, there are three minors
of
order $m$
(up to signs):
$$
\det\pmatrix{
4&5\cr
7&5\cr
}=-15,\quad
\det\pmatrix{
4&5\cr
-2&0\cr
}=10,\quad
\det\pmatrix{
7&5\cr
-2&0\cr
}=10,
$$
and six minors of order $m-1$: 4, 5, 7, 5, $-2$, 0. Thus,
the theorem asserts
that (in this example) the number of solutions is divisible by
$$
\GCD\({\GCD(-15,10,10)\over\GCD(4, 5, 7, 5, -2, 0)},
|C(a)\cap C(b)|\)=
\GCD(5,|C(a)\cap C(b)|).
$$

Note that the agreements about boundary cases in Theorem 1 are natural.
Indeed, we always can add a fictitious equation 1=1 to make the number
of equations larger than $m$.  We can also add a new variable $z$ and the
equation~$z=1$ (this does not affect the number of solutions and makes
$m>1$).  As for the philosophical question on the interpretation of the
fraction $0\over0$, it can be understood arbitrarily, e.g., the reader may
assume that ${0\over0}=\the\year$; in any case, Theorem 1 remains valid
(but weaker than under the suggested interpretation).

The meaning of the value~$\Delta_m\over\Delta_{m-1}$ is as follows.  It is
well known (see, e.g., [Vin03]) that invertible integer elementary
transformations of rows and columns can transform any integer matrix $A$
into a diagonal matrix, where the diagonal entries divide each other
(each diagonal entry divides the next one). This diagonal matrix is
uniquely determined up to the signs of diagonal elements (and is sometimes
called the \emph{Smith form} of $A$); the diagonal elements of the Smith
form (sometimes called the \emph{invariant factors} of $A$) equal to the
ratios ${\Delta_i\over\Delta_{i-1}}$. Thus, in these terms,
$\Delta_m\over\Delta_{m-1}$ is the 
$m$-th 
invariant factor of the matrix of
the system of equations.  One can also say that
\disp{\sl
the absolute
value of $\Delta_m\over\Delta_{m-1}$ is the period \(exponent\) of the
quotient of the free abelian group $\Z^m$ by the subgroup generated by the
rows of the matrix of the system of equations}
(with the stipulation that
this ratio vanishes if and only if the period is infinite).

The Frobenius and Solomon theorems as well as their generalisations stated
above are special cases of Theorem~1.

The following theorem is on the first view similar to neither the
Frobenius theorem nor the Solomon theorem.

\proclaim
Iwasaki theorem \rm [Iwa82].
For any integer $n$, the number of elements of a finite group $G$ whose
$n$-th 
powers lie in a subgroup $H\subseteq G$ is divisible by $|H|$.

This beautiful theorem remains (for some reason) not widely known. In
[SaAs07], it was noticed that the divisibility by~$|H|$ still holds for
the number of solutions to the ``equation" $x^n\in HgH$, where $HgH$ is
any double coset of a subgroup~$H$. Clearly, the Iwasaki theorem and its
generalisations deals with predicates that are not equations in the usual
sense.  Let us say that a \emph{generalised equation} over a group $G$ is
an expression of the form~$w(x_1,\dots,x_n)\in HgH$, where $H$ is a
subgroup of~$G\ni g$, and $w(x_1,\dots,x_m)$ is an element of the free
product $G*F(x_1,\dots,x_m)$ of $G$ and a free group; in other terms, $w$
is a word in the alphabet~$G\sqcup\{x_1^{\pm1},\dots,x_m^{\pm1}\}$. The
elements of $G$ 
occurring 
in this word are called the \emph{coefficients}
of the generalised equation.  A system of generalised equations, a
solution to this system, and a matrix of this system are defined in a
natural way.

In [KM17], the following generalisation of the Iwasaki theorem was obtained:
\disp{%
\hfuzz500mm
\sl
the number of solutions to a system of generalised coefficient-free
equations whose right-hand sides are double cosets of the same
subgroup~$H$ {\rm(e.g.,
$
\{x^{100}y^{\the\year}[x,y]^4\in Hg_1H,
\quad
  [x^5,y^6]^7(xy)^8\in Hg_2H,\dots
\}
$)}
is divisible by $|H|$.
}

The following theorem includes all results stated above.

\Th 2%
\fn{%
{\bf Theorem 1} in the journal version.}%
.
Let $S$ be a
\(not necessarily finite\)
system of
generalised equations
in finitely many unknowns $x_1,\dots,x_m$
over a group $G$
and let $P$ be its subsystem:
$$
S=\{u_i(x_1,\dots,x_m)\in H_ig_iH_i\;|\;i\in I\}
\supseteq
P=\{u_j(x_1,\dots,x_m)\in H_jg_jH_j\;|\;j\in J\},
$$
\(where $J\subseteq I$,\quad
$u_i\in G*F(x_1,\dots,x_m)$,\quad
$g_i\in G$, and $H_i$ are subgroups of $G$\).
Then the number of solutions to~$S$
in~$G$
is divisible by the
greatest common divisor of the
subgroup
$$
\~H=
\left(\bigcap\limits_{j\in J}N(H_jg_jH_j)\right)
\cap
\left(\bigcap\limits_{i\in I\setminus J}H_i\right)
\cap
\bigl(\hbox{\rm the centraliser of the set of
                coefficients of $S$}\bigr)
$$
and the number~$\Delta_m\over\Delta_{m-1}$,
where $\Delta_k$ is the greatest common divisor
of
all minors of order $k$ of the matrix of the subsystem~$P$.
\rm Henceforth, $N(A)\:=\{g\in G\;|\;g^{-1}Ag=A\}$ is the normaliser of
a subset $A$ in a group~$G$.

To deduce Theorem 1 from Theorem 2,
we
rewrite the system of equations in the ``generalised" form, i.e.
we put
$
S=P=\Bigl\{u_1(x_1,\dots,x_m)\in\11\1,
\quad
u_2(x_1,\dots,x_m)\in\11\1,\dots\Bigr\}
$
and note that the normaliser of the trivial
subgroup is the whole group.

On the other hand, setting
$$
S=\Bigl\{u_1(x_1,\dots,x_m)\in Hg_1H,
\quad
u_2(x_1,\dots,x_m)\in Hg_2H,\dots\Bigr\}
\qqbox{and}
P=\emptyset
\qbox{(where $u_i\in F(x_1,\dots,x_m)$)},
$$
we obtain
the mentioned above generalisation (from [KM17]) of the Iwasaki theorem.


As a matter of fact, a relation between Solomon's and Iwasaki's theorems
was established in [KM14] and [KM17]; our achievement consists only in
adding ``Frobeniusness".  The main theorem of [KM17] says that, if we have
a group~$F$ with a fixed epimorphism onto~$\Z$ and some set of
homomorphisms from $F$ into another group $G$, and this set is invariant
with respect to some natural transformations (depending on the epimorphism
$F\to\Z$ and a subgroup $H$ of $G$), then the number of these
homomorphisms $F\to G$ is divisible by $|H|$. Choosing suitable
sets of homomorphisms, the authors of [KM17] obtained Solomon's and
Iwasaki's theorem as special cases of their main theorem.

Our main theorem (see Section 1) is a modular analogue of the main theorem
of [KM17]:  we take an epimorphism~$F\to\Z/n\Z$ instead of $F\to\Z$.  One
can say that the main theorem of this paper is related to the main theorem
of [KM17] in the same way as Theorem 1 to the generalisation (from [KM14])
of the Solomon theorem mentioned in the beginning of this paper. An
important role in our argument is played by an elementary (but nontrivial)
lemma due to Brauer [Bra69]. Actually, we need this lemma not to prove the
main theorem but rather to explain that its statement \emph{per se} makes
some sense.  For readers' convenience, we give a proof of the Brauer lemma
in the last section. Section~5 contains the proof of the main theorem.

In Section~2, we deduce Theorem 2 from the main theorem. As another
corollary, we obtain a theorem on  equations in rings (Theorem~3 in
Section~3) that implies, e.g., the following fact, which can be considered
as a generalisation of the Frobenius theorem in another
direction:

{\sl
for any representation $\rho\:G\to\GL(V)$ of a
group $G$ and any words
$u_i(x_1,\dots,x_m)
\in F(x_1,\dots,x_m)$,
$$
\hfuzz8.5pt
\hbox{the number of solutions to the equation }
\sum_{i=1}^k\Bigl(\rho\bigl(u_i(x_1,\dots,x_m)\bigr)\Bigr)^{l_i}
=\id
\hbox{ is divisible by}
\cases{
\!\!\GCD\bigl(G,\GCD(\{l_i\})\bigr) & \kern-7pt  always;       \cr
\!\!\GCD\bigl(G,\LCM(\{l_i\})\bigr) & \kern-7pt if $k\le m$;   \cr
\!\!|G|                             & \kern-7pt if $k<m$.      \cr
}
$$
}

In Section 4, we show that the main theorem
implies
some
fact about the number of crossed homomorphisms, generalising earlier known
results.
In the next to last section, we discuss open questions.


The authors thank Saveliy Skresanov for valuable remarks.


\goodbreak

{\noindent \bf Notation and conventions}
we use are mainly standard. Note only that, if
$k\in \Z$ and $x$ and $y$ are elements of a group, then $x^y$,
$x^{ky}$, and $x^{-y}$ denote
$y^{-1}xy$, $y^{-1}x^ky$ and $y^{-1}x^{-1}y$,
respectively.
The commutator subgroup of a group~$G$ is denoted by $G'$ or $[G,G]$.  If
$X$ is a subset of a group, then $|X|$, $\gp X$, $\nc X$, $C(X)$,
and~$N(X)$ are the cardinality of~$X$, subgroup generated by~$X$, normal
closure of~$X$, centraliser of~$X$, and normaliser of~$X$. The index of a
subgroup $H$ of a group $G$ is denoted by $|G:H|$. The letter~$\Z$ denotes
the set of integers.  If $R$ is an associative ring with unity, then $R^*$
denotes the group of units of this ring. $\GCD$ and $\LCM$ are the
greatest common divisor and least common multiple.  The symbol $\exp(G)$
denotes the period (exponent) of a group $G$ if this period is finite; we
assume $\exp(G)=0$ if the period is infinite.  The symbol~$\gp g_n$
denotes the cyclic group of order $n$ generated by an element $g$. The
free group of rank~$n$ is denoted by $F(x_1,\dots,x_n)$ or $F_n$.  The
symbol $A*B$ denotes the free product of groups $A$ and~$B$.

Let us recall once again that
the finiteness of groups is not
assumed by default; the divisibility is always understood in the sense
of cardinal arithmetics (an infinite cardinal is divisible by all nonzero
cardinals not exceeding it), and
$
\GCD(G,n)\:=\LCM\(\bigl\{|H|\;\bigm|\;
\hbox{$H$ is a subgroup of $G$, and $|H|$ divides $n$}\bigr\}\).
$

\s 1.
Main theorem

A group $F$ equipped with an epimorphism $F\to\Z/n\Z$ (where $n\in\Z$) is
called an \emph{$n$-indexed} group. This epimorphism~$F\to\Z/n\Z$ is
called \emph{degree} and denoted $\deg$. Thus, to any element~$f$ of an
indexed group $F$, an element~$\deg f\in\Z/n\Z$ is assigned; the group $F$
contains elements of all degrees and $\deg(fg)=\deg f+\deg g$ for
any~$f,g\in F$.

Suppose that $\phi\:F\to G$ is a homomorphism from an $n$-indexed group $F$
to a group $G$ and $H$ is a subgroup of~$G$.
The subgroup
$$
H_\phi=\bigcap_{f\in F}H^{\phi(f)}\cap C(\phi(\ker\deg))
$$
is called the
\emph{$\phi$-core} of $H$ [KM17]. In other words, the
$\phi$-core $H_\phi$ of $H$ consists of
elements~$h$ such that $h^{\phi(f)}\in H$ for all $f$, and
$h^{\phi(f)}=h$ if $\deg f=0$.

\proclaim{Main theorem}.
Suppose that an
integer $n$
is a multiple of the order of a subgroup
$H$ of group $G$ and a
set $\Phi$ of homomorphisms from an $n$-indexed group $F$ to $G$
satisfies the following
conditions.
\item{\rm I.}
$\Phi$ is invariant with respect to conjugation by elements of $H$:
$$
\qbox{if $h\in H$ and $\phi\in\Phi$, then the homomorphism
$\psi\:f\mapsto\phi(f)^h$ lies in $\Phi$.}
$$
\item{\rm II.}
For any $\phi\in\Phi$ and any element $h$ of the $\phi$-core
$H_\phi$ of $H$, the homomorphism $\psi$ defined by
$$
\psi(f)=
\cases{
\phi(f)& for all elements $f\in F$ of degree zero;
\cr
\phi(f)h& for
some element
$f\in F$ of degree one
\small(and, hence, for all degree-one elements)
\cr
}
$$
belongs to $\Phi$ too.
\enditem
Then $|\Phi|$ is divisible by
$|H|$.

Note that the mapping $\psi$ from Condition I is a homomorphism for any
$h\in G$, and the formula for $\psi$ from Condition~II defines a
homomorphism for any $h\in H_\phi$ (as explained below). Thus,
Conditions~I and II only require these homomorphisms to belong to $\Phi$.

\Lemma 0%
\fn{%
{\bf Lemma 2} in the journal version.}%
.
Suppose that $\phi\:F\to G$ is a
homomorphism from an
$n$-indexed group $F$ to a group $G$,
$f_1\in F$ is an element of degree one
and $g\in G$. Then the
homomorphism $\psi\:F\to G$
such that $\psi(f)=\phi(f)$ for all $f\in F$ of degree zero and
$\psi(f_1)=\phi(f_1)g$ exists if and only if
$g\in C(\phi(\ker\deg))$ and
$\bigl(\phi(f_1)g\bigr)^n=\bigl(\phi(f_1)\bigr)^n$.

\Proof
The group $F$ can be presented in the form
$$
F\iso\bigl(F_0*\gp{x}_\infty\bigr)/
\nc{\{u^xu^{-f_1}\;|\; u\in F_0\}\cup\{x^nf_1^{-n}\}},
\qbox{where } F_0=\ker\deg.
$$
Therefore,
the mapping $\psi\:F_0\cup\{x\}\to G$ can be extended to a homomorphism
if and only if
its restriction to $F_0$ is homomorphism and
the relations
$u^x=u^{f_1}$ (for $u\in F_0$) and $x^n=f_1^n$ are mapped to true
equalities in $G$:
$$
\psi(u)^{\psi(x)}=\psi(u^{f_1})
\qqbox{and}
\psi(x)^n=\psi(f_1^n).
\eqno{(*)}
$$
If the restrictions of $\psi$ and $\phi$ to $F_0$ coincide and
$\psi(x)=\phi(f_1)g$, then the first equality $(*)$ says that $g$ commutes
with~$\phi(u)$ (for all $u\in F_0$), while the second equalities $(*)$
takes the form $\bigl(\phi(f_1)g\bigr)^n=\bigl(\phi(f_1)\bigr)^n$.  This
completes the proof.

\smallskip

Recall also the following beautiful (but not widely known) fact.

\proclaim
Brauer lemma \rm[Bra69].
If $U$ is a finite normal subgroup of a group $V$, then, for all $v\in V$
and $u\in U$, the elements~$v^{|U|}$ and $(vu)^{|U|}$ are conjugate by an
element of $U$.

\smallskip

These two lemmata imply immediately that the mapping $\psi$ from
Condition~II is a homomorphism
for any $h\in H_\phi$ because
$(\phi(f)h)^n=(\phi(f))^n$ by the Brauer lemma applied to
$
U=H_\phi\subset V=H_\phi\cdot\gp{\phi(f_1)}\ni \phi(f_1)=v.
$
Indeed, we obtain the equality
$\bigl(\phi(f_1)h\bigr)^{|H_\phi|}=\bigl(\phi(f_1)\bigr)^{|H_\phi|u}$
for some $u\in H_\phi$ and, hence,
$
\bigl(\phi(f_1)h\bigr)^n=\bigl(\phi(f_1)\bigr)^{nu}=
\bigl(\phi(f_1^n)\bigr)^u
$
(because
$|H_\phi|$ divides $n$).  It remains to note that $u\in H_\phi$ commutes
with $\phi(f_1^n)$ because $\deg f_1^n=n=0\in\Z/n\Z$. Thus, we obtain the
equality $\bigl(\phi(f_1)h\bigr)^n=\bigl(\phi(f_1)\bigr)^n$.  It remains
to refer to Lemma 0.

\medskip

In the case $n=0$ the main theorem was proved in [KM17]. So, our
theorem is a ``modular analogue" of the main result of [KM17].  On the
other hand, our main theorem is deduced (in Section~5) from this special
case $n=0$.

\Lemma 1%
\fn{%
{\bf Lemma 3} in the journal version.}%
.
In Condition II of the main theorem, $\psi(f)\in\phi(f)H_\phi$ for all
$f\in F$.

\Proof
Indeed, if $\deg f=d$, then $f=f_1^df_0$,
where $f_1$ is the (fixed) element of degree one
(from Condition II) and $f_0$ is an element of degree zero. Then
$$
\psi(f)=\psi(f_1)^d\psi(f_0)=(\phi(f_1)h)^d\phi(f_0)
=\!=\!=
\phi(f_1)^d\phi(f_0)h'=\phi(f_1^df_0)h'=\phi(f)h',
$$
where the equality $=\!=\!=$ is valid for some $h'\in H_\phi$
because $h\in H_\phi$ and $\phi(F)$ normalises $H_\phi$.

\s 2.
Proof of Theorem 2

Let $L\subseteq G$ by the subgroup generated by all
coefficients of the system $S$.
Take as $H$ any subgroup of the group $\~H$
whose
order divides $n\:={\Delta_m\over\Delta_{m-1}}$,
and put
$$
F=L*F(x_1,\dots,x_m)
\qqbox{and}
\Phi=\Bigl\{\phi\:F\to G\;\Bigm|\;
\phi(f)=f \hbox{ for $f\in L$}
\qqbox{and}
\phi(u_i)\in H_ig_iH_i
\hbox{ for $i\in I$}\Bigr\}.
$$
As the indexing $\deg\:F\to\Z/n\Z$, take an epimorphism
whose kernel contains
$L$ and all~$u_j$, where $j\in J$. Such an epimorphism exists
because $n$ is the period of the finitely generated abelian
group~$F/\bigl([F,F]\cdot L\cdot\gp{\{u_j\;|\;j\in J\}}\bigr)$.

Let us verify that the conditions of the main theorem hold.
Condition I holds obviously for all $h\in H$ (and even for all
$h\in\~H$) because (by definition) $\~H$  centralises
$L$ and normalises double cosets $H_ig_iH_i$.

\noindent
Condition II holds also for all $h\in H_\phi$ because
\-
on $L$, the homomorphism $\psi$ coincides with $\phi$
as $L$ consists of zero-degree elements;
\-
$\psi(u_j)=\phi(u_j)$ for $j\in J$ because again
$\deg u_j=0$;
\-
for $i\in I\setminus J$, we have
$\psi(u_i)\in \phi(u_i)H_\phi\subseteq \phi(u_i)H_i$
(where the inclusion $\in$ follows from Lemma 1).

\enditem
Thus, the main theorem implies that $|\Phi|$ is divisible by the order of
any subgroup $H\subseteq\~H$ whose order divides~$n$, i.e. $|\Phi|$ is
divisible by $\GCD(\~H,n)$. It remains to note that $|\Phi|$ is the
number of solutions to $S$.

\s 3.
Rings and representations

A \emph{generalised homogeneous modulo $n$}
equation with a set of unknowns~$X$ over an associative unital ring
$R$ is a finite expression of the form
$$
\sum_i\prod_j c_{ij}x_{ij}^{k_{ij}}=0,
\qbox{where
\emph{coefficients} $c_{ij}\in R$,
\emph{unknowns} $x_{ij}\in X$,
and \emph{exponents} $k_{ij}\in\Z$},
$$
such that, for some mapping $\deg\:X\to\Z/n\Z$, the value
$\sum\limits_jk_{ij}\deg(x_{ij})$ (called the \emph{degree of
the equation}) does not depend on $i$ (i.e.  the ``polynomial" in the
left-hand side of the equation is homogeneous with respect to some
assigning of degrees to variables), and
$\gp{\{\deg x\;|\;x\in X\}}=\Z/n\Z$.

A system of equations is called generalised homogeneous modulo $n$ if all
equations of this system are generalised homogeneous modulo $n$
(of possibly different degrees) with respect to the same
function~$\deg\:X\to\Z/n\Z$.

As we explain below, the set
$
M=\{n\in\Z\;|\;
\hbox{a given system is generalised homogeneous modulo $n$}\}
$
consists of all divisors of a number $n_0$, called
the \emph{homogeneity modulus} of the system. In other words,
the homogeneity modulus
is the maximal number from $M$ or zero if $M$ is infinite.

To find the homogeneity modulus, consider a
homogenous
system of linear
equations, where unknowns are degrees of variables
and also (the negations of) degrees of equations;
these linear equations say that the degree of each monomial
equals the degree of the corresponding equation.
The matrix of this
system (called the \emph{homogeneity matrix}
of the initial system of equations) has the following form.
Suppose that $X=\{x_1,\dots,x_m\}$.
The \emph{homogeneity matrix of 
the $p$-th
equation} is the integer matrix $A_p=(a_{kl})$
of size
$$
(\hbox{the total number of monomials in the system})
\times
\bigl(
m
+
(\hbox{the number of equations})
\bigr),
$$
where, for $l\le m$, the
$(k,l)$th entry is the exponent sum of 
the $l$-th 
unknown in
the $k$-th 
monomial, the 
$(m+p)$-th 
column consists of ones, and the remaining
columns are zero for $l>m$. The homogeneity matrix of the system of
equations is composed from the matrices $A_p$ written one under another:
$
A=\pmatrix{A_1\cr A_2\cr \vdots}.
$
For example, the system of equations
$
\left\{
ax^3y^2+y^7bx-1=0,
\quad
xy^2x+y^7x^5=0
\right\}
$
(where $x$ and $y$ are
unknowns and $a,b\in R$ are coefficients)
has the following homogeneity matrix:
$$
A=
\pmatrix{
3 & 2 & 1 & 0\cr
1 & 7 & 1 & 0\cr
0 & 0 & 1 & 0\cr
2 & 2 & 0 & 1\cr
5 & 7 & 0 & 1\cr
},
\qqbox{composed of matrices}
A_1=
\pmatrix{
3 & 2 & 1 & 0\cr
1 & 7 & 1 & 0\cr
0 & 0 & 1 & 0\cr
}
\hbox{ and }
A_2=
\pmatrix{
2 & 2 & 0 & 1\cr
5 & 7 & 0 & 1\cr
}.
$$

\proclaim
Homogeneity-modulus lemma.
The homogeneity modulus of a system of
$s$ equations
in $m$ unknowns
over
an associative ring with unity
is
$\frac{\Delta_{m+s}}{\Delta_{m+s-1}}$, where~$\Delta_i$ is the greatest
common divisor of all minors of order $i$ of the homogeneity matrix of the
system. As always, the following
conventions are assumed: $\Delta_i=0$ if the total number of monomials in
all equations is less than $i$; $\Delta_{0}=1$; $\frac{0}{0}=0$.

\Proof
Let $A$ be the homogeneity matrix. We have to find the maximal
number $n$ such that the system of linear homogeneous equations $AX=0$
(in $m+s$ variables) has a solution in $\Z/n\Z$ whose components
generate $\Z/n\Z$ as an additive group
(or, equivalently, the first $m$ components of the solution
generate $\Z/n\Z$, because
the equations say that the last $s$ components are combinations of the
first $m$ ones).
In other words, $n$ is the largest order of cyclic quotient of the
finitely generated group $\Z^{m+s}/N$, where
$N$ is the subgroup generated by rows of $A$.
As noted already,
the largest cyclic quotient $n$ of $\Z^{m+s}/N$ is
$\frac{\Delta_{m+s}}{\Delta_{m+s-1}}$, as required.

\Th 3%
\fn{%
{\bf Theorem 4} in the journal version.}%
.
Let $R$ be an associative ring with unity and let $G$ be a subgroup
of the
multiplicative group 
of
this ring. Then, for each system of equations
over $R$ in $m$ unknowns, the number of its solutions lying in $G^m$ is
divisible by the greatest common divisor of the homogeneity modulus
of the system and the intersection of $G$ with the centraliser of the set
of coefficients of the system.

\Proof
Let $G_0$ be the intersection of $G$ and the centraliser of the set
of coefficients and let $n$ be the homogeneity modulus.
Consider the free group $F=F(X)$ (where $X$ is the set of unknowns)
and an epimorphism $\deg\:F\to\Z/n\Z$.

Let us apply the main theorem taking
$\Phi$ to be the set of all
homomorphisms $\phi\:F\to G$ such that the
tuple~$(\phi(x_1),\dots,\phi(x_m))$ is a solution to the system of
equations (so, the number of solutions is $|\Phi|$). Take $H$ to be any
subgroup of $G_0$ of order dividing~$n$. Condition I of the main theorem
obviously holds. To verify Condition II, choose an element $t\in F$ of
degree one and write each variable $x_i$ in the form $x_i=t^{\deg x_i}y_i$,
where $y_i=t^{-\deg x_i}x_i$ has degree zero.  In new notation, each
equation $w(x_1,\dots,x_m)=0$ takes the form
$v(t,y_1,\dots,y_m)=0$ and
the exponent sum of $t$
in each term of
this equation
is the same (modulo $n$).
Now, note that, if
$v(\phi(t),\phi(y_1),\dots,\phi(y_m))=0$
and~$h\in H_\phi$, then
$v(\phi(t)h,\phi(y_1),\dots,\phi(y_m))=0$. This follows from the
(right) divisibility of $v(\phi(t)h,\phi(y_1),\dots,\phi(y_m))$
by~$v(\phi(t),\phi(y_1),\dots,\phi(y_m))$ due to the following fact.

\proclaim
Fact {\rm([KM17], Lemma 1)}.
If $M$ is a monoid, $b_i,a,h\in M$, elements $a$ and $h$ are invertible,
and the
elements $a^{-s}ha^s$, where~$s\in\Z$, commute with all $b_i$, then,
for any expression of the form
$
u(t)=b_0t^{m_1}b_1\dots t^{m_l}b_l,
\qbox{where $m_i\in\Z$,}
$
we have
$
u(ah)=
\cases{
h^{a^{-1}}h^{a^{-2}}\dots h^{a^{-k}}u(a) & if $k=\sum m_i>0$;
\cr
h^{-1}h^{-a}\dots h^{-a^{-1-k}}u(a) & if $k=\sum m_i<0$;
\cr
u(a), & if $k=\sum m_i=0$.
\cr
}
$

We apply this fact to each term of $v$;
we also use that $t^n$ has
degree zero and $(\phi(t)h)^n=(\phi(t))^n$ according to Lemma 0.

Thus, the main theorem implies that $|\Phi|$ (i.e. the number of solutions
to the system of equations) is divisible by~$|H|$ as required
(because $H$ is an arbitrary subgroup of $G_0$ whose order divides
the homogeneity modulus).

\Example.
{\sl
If $\rho\:G\to R^*$ is a homomorphism from a finite group $G$ to
the multiplicative group of an associative ring $R$ with unity
{\rm(e.g.,
$\rho\:G\to\GL(V)$ is a linear representation of $G$)}, then,
for
any words
$u_i(x_1,\dots,x_m)
\in F(x_1,\dots,x_m)$,
$$
\hfuzz8.5pt
\hbox{the number of solutions to the equation }
\sum_{i=1}^k\Bigl(\rho\bigl(u_i(x_1,\dots,x_m)\bigr)\Bigr)^{l_i}
=1
\hbox{ is divisible by}
\cases{
\!\!\GCD\bigl(G,\GCD(\{l_i\})\bigr) & \kern-7pt  always;       \cr
\!\!\GCD\bigl(G,\LCM(\{l_i\})\bigr) & \kern-7pt if $k\le m$;   \cr
\!\!|G|                             & \kern-7pt if $k<m$.      \cr
}
$$
}
To show this, it suffices to
apply Theorem 3 to the subgroup
$\rho(G)\subseteq R^*$.
The homogeneity matrix of this equation has the form
$
B=\pmatrix{
A&\!\!\!\!\!\matrix{1\cr\vdots\cr}\cr
0\ \dots\ 0
&\!\!\!\!\!1\cr
},
$
where the last row corresponds to 1 in the right-hand side of the equation,
and the 
$i$-th 
row of the matrix $A$ corresponds to the 
$i$-th 
term in the
left-hand side of the equation and, therefore, all elements of this row
are divisible by $l_i$. It remains to note that the 
$j$-th 
invariant
factor of the matrix $B$ coincides with the 
$(j-1)$-th 
invariant factor
of $A$ and use the following fact, which we leave to readers as an easy
exercise:  

{\sl\noindent
if the 
$i$-th 
row of an integer matrix $k\times m$ is divisible by $l_i$,
\newline
\phantom{1}\quad\quad\quad\quad\quad\quad\quad\quad
then
the 
$m$-th 
invariant factor of this matrix%
$\cases{
\hbox{is divisible by }\GCD(\{l_i\})& always;\cr
\hbox{is divisible by }\LCM(\{l_i\})& for $k=m$;\cr
\hbox{vanishes}& for $k<m$.\cr
}
$
}

\medskip

Note that Theorem 1 can be obtained as a corollary of Theorem 3.
Indeed, take $R=\Z G$; the group ring
contains $G$ as a subgroup of the
multiplicative group. Any system of equations over $G$ can be rewritten in
``ring"\ form: $\{w_i(x_1,\dots)-1=0\}$. It remains to note
that the value $\Delta_m\over\Delta_{m-1}$ from Theorem 1 becomes
exactly the homogeneity modulus from the homogeneity-modulus lemma.

\s 4.
Crossed homomorphisms

Suppose that a group $F$ acts (on the right) on a group $B$
by automorphisms: $(f,b)\mapsto b^f$. Recall that a
\emph{crossed homomorphism} from $F$ to $B$ with respect to this action is
a mapping $\alpha\:F\to B$ such that
$\alpha(ff')=\alpha(f)^{f'}\alpha(f')$ for all~$f,f'\in F$.
Saveliy Skresanov noted that the main theorem easily implies
the following fact proved in [ACNT13]
(using character theory)
for finite groups
$F$ and $B$.

\Th 4%
\fn{%
{\bf Theorem 5} in the journal version.}%
.
If a group $F$ admitting an epimorphism onto $\Z/n\Z$ acts by
automorphisms on a group $B$, then the number of crossed homomorphisms
$F\to B$ is divisible by $\GCD(B,n)$.

\Proof
The set of crossed homomorphisms is in one-to-one correspondence with the
set $\Phi$ of (usual) homomorphisms from $F$ to the semidirect
product~$G=F\semitimes B$ (with respect to the given action) such that
their compositions with the projection~$\pi\:F\semitimes B\to F$ is the
identity mapping~$F\to F$. We have to show that $|\Phi|$ is a multiple of
$|H|$ for any subgroup $H\subseteq B$ whose order divides $n$ (by
definition of $\GCD(B,n)$).

The group $F$ is $n$-indexed by the hypothesis of Theorem 4. Therefore,
the assertion follows immediately from the main theorem.  Conditions of
the main theorem hold by trivial reasons:  Condition I is fulfilled
because $\pi(h^{-1}gh)=\pi(g)$; Condition II follows immediately from
Lemma 1 because~$\pi(gh)=\pi(g)$ (for $g\in G$ and $h\in H$).

\s 5.
Proof of the main theorem

Take an element $f_1\in F$ of degree one,
put $F_0=\ker\deg\subset F$, and consider
the semidirect product
$\~F=\gp a_\infty\semitimes F_0$, where $a$ acts on
$F_0$ as $f_1$ does:
$
u^a=u^{f_1} \hbox{ for $u\in F_0$}.
$
The group
$\~F$ admits a natural indexing (0-indexing) $\deg\:\~F\to\Z$ (denoted
by the same symbol~$\deg$). The kernel of this map is $F_0$ and
$\deg a=1$. Moreover, there is a natural epimorphism~$\alpha\:\~F\to F$
mapping $a$ to $f_1$ and identity on~$F_0$. Let us verify that
the conditions of the main theorem hold for the set
$\~\Phi=\{\phi\circ\alpha\;|\;\phi\in\Phi\}$
of homomorphisms from $\~F$ to $G$.

Condition I holds obviously. To verify Condition II,
take the degree-one element~$a\in \~F$ and some
homomorphism $\~\phi=\phi\circ\alpha\in\~\Phi$
(where $\phi\in\Phi$).
Then
the homomorphism $\~\psi$ from
Condition~II has the form
$$
\~\psi(\~f)=
\cases{
\phi(\~f)& for all elements $\~f\in F_0$;
\cr
\phi(f_1)h& for $\~f=a$;
\cr
}
\quad\qbox{where $\phi\in\Phi$ and $h\in H_{\tphi}$.}
\eqno{(1)}
$$
We have to show that
$\~\psi$ lies in $\~\Phi$,
i.e. has the form $\~\psi=\phi'\circ\alpha$, where $\phi'\in\Phi$.
Note that
$
H_\tphi=H_\phi,
$
because the images of $\~\phi=\phi\circ\alpha$ and $\phi$
coincide, and the images of
zero-degree elements for these homomorphisms
coincide:
$
\~\phi(\ker\deg)=\~\phi(F_0)=\phi(F_0).
$
Formula (1) takes the form
$$
\~\psi(\~f)=
\cases{
\phi(\~f)& for $\~f\in F_0$;
\cr
\phi(f_1)h& for $\~f=a$;
\cr
}
\quad\qbox{where $\phi\in\Phi$ and $h\in H_\phi$.}
$$
This means that $\~\psi=\psi\circ\alpha$, where
$$
\psi(f)=
\cases{
\phi(f)& for $f\in F_0$;
\cr
\phi(f_1)h& for $f=f_1$;
\cr
}
\quad\qbox{where $\phi\in\Phi$ and $h\in H_\phi$.}
$$
The homomorphism $\psi\:F\to G$ lies in $\Phi$ by Condition II of the
theorem we are proving. Therefore, $\~\psi\in\~\Phi$. Thus, the conditions
of the
main theorem hold
for the set $\~\Phi$ of homomorphisms from the 0-indexed
group~$\~F$ to $G$.
Therefore,
$|\~\Phi|$ is divisible on $|H|$
by virtue of the main theorem of [KM17].
It remains to note that $|\Phi|=|\~\Phi|$
since $\alpha$ is surjective. This completes the proof.


Note that we do not verify here that $\psi$
defines a homomorphism; this is non-obvious but true,
see Section 1.


\s 6.
Open questions

Theorems 1,2,3,4 assert that some numbers are multiples of the ratios of
two integers. Oddly, we do not know {\sl whether these ratios can be
replaced by their numerators.}

\Questions 1 and 2%
\fn{%
{\bf Questions 6 and 7} in the journal version.}%
.
Is it possible to replace the ratio
$\Delta_m/\Delta_{m-1}$ by its numerator $\Delta_m$ in Theorems 1 and 2?

For coefficient-free systems of equations,
Question 1 is equivalent to
the following question
posed in [AsYo93] (for finite groups $F$ and $G$):
\disp{\sl
\hfuzz11cm
is the number of homomorphisms from a finitely
generated
group $F$ to a
group $G$ divisible by~$\GCD(|F/F'|,G)$\?
}%
This problem remains unsolved even for finite groups (as far as we know).
A survey of some results can be found in~[AsTa01]; e.g., the answer is
positive if $F$ is abelian [Yosh93].

Theorem 3 suggests a similar question.

\Question 3%
\fnn{%
{\bf Question 8} in the journal version.}%
.
Is it possible, in Theorem 3, to replace
the homogeneity modulus by
its numerator $\Delta_{m+s}$
\rm(see the homogeneity-modulus lemma)?

As for Theorem 4,
it also leads us to
a similar question.
Indeed, Theorem 4 implies, in particular,
that
{\sl
if a finitely generated group $F$
acts by automorphisms on a group $B$,
then the number
of crossed homomorphisms $F\to B$ is divisible by $\GCD(\exp(F/F'),B)$.
}

\Question 4%
\fnnn{%
{\bf Question 9} in the journal version.}%
.
Is it possible, in the proposition above, to replace
the period $\exp(F/F')$ by the order
of this quotient group?

This question was posed
for the first time in [AsYo93]
(for finite groups $F$ and $B$).
To show the similarity of Questions 4 and~1, we
recall that the absolute value of the ratio $\Delta_m/\Delta_{m-1}$ in
Question~1 is the period of the quotient group of the free abelian group
$\Z^m$ by the subgroup generated by the rows of the matrix of the system
of equations, while the absolute value of the numerator $\Delta_m$ is the
order of this quotient group.

\s 7.
Proof of the Brauer lemma

We follow the original proof from
[Bra69] but translate it into a more convenient (in our view) language.

\proclaim
Brauer Lemma \rm[Bra69].
If $U$ is a finite
normal subgroup of a group $V$,
then, for all $v\in V$ and $u\in U$, the
elements~$v^{|U|}$~and~$(vu)^{|U|}$ are conjugate
by an element of $U$.

\Proof
The group $\Z$ acts by permutations on the subgroup $U$:
$$
a\circ i=v^{-i}a(vu)^i,
\qbox{(where $i\in\Z \hbox{ and } a\in U$)}.
$$
Let $m$ be the minimum length of an orbit.
In other words, $m$ is the minimum length of a cycle
in the decomposition of the permutation $a\mapsto v^{-1}avu$
(of $U$) into the product of independent cycles.
The set
$
X=\{a\in U\;|\; a\circ m=a\}
$
is the union of all orbits of length $m$;
therefore, $|X|$ is divisible by $m$.
On the other hand, (by definition of the action)
$
X=\{a\in U\;|\; v^{-m}a(vu)^m=a\}=\{a\in U\;|\; a^{-1}v^ma=(vu)^m\}
$
and, hence, $|X|$ is the order of the centraliser of $v^m$ in~$U$
(because,
in any group,
a nonempty
set of the form $\{x\;|\; x^{-1}yx=z\}$ is a coset of the
centraliser of $y$).
Thus, $|X|$ divides $|U|$ and, therefore, $m$ divides
$|U|$ and $a\circ|U|=a$ (if $a$ lies in an orbit of length $m$).
This competes the proof.


\References



[Vin03]
E. B. Vinberg,
A Course in Algebra
(Issue 56 of Graduate studies in mathematics),
American Mathematical Soc., 2003.


[Stru95]
S. P. Strunkov,
On the theory of equations in finite groups,
Izvestiya: Math., 59:6 (1995), 1273-1282.


[AmV11]
A. Amit, U. Vishne,
Characters and solutions to equations in finite groups,
J. Algebra Appl., 10:4 (2011), 675-686.

[And16]
R. Andreev,
A translation of
``Verallgemeinerung des Sylow'schen Satzes"
by F. G. Frobenius.
\newline
arXiv:1608.08813.

[ACNT13]
T. Asai, N. Chigira, T. Niwasaki, Yu. Takegahara,
On a theorem of P. Hall,
Journal of Group Theory, 16:1 (2013), 69-80.

[AsTa01]
T. Asai, Yu. Takegahara,
$|\Hom(A,G)|$, IV,
J. Algebra, 246 (2001), 543-563.

[AsYo93]
T. Asai, T. Yoshida,
$|\Hom(A,G)|$, II,
J. Algebra, 160 (1993), 273-285.

[Bra69]
R. Brauer,
On A Theorem of Frobenius,
The American Mathematical Monthly, 76:1 (1969), 12-15.


[BrTh88]
K. Brown, J. Th\'evenaz,
A generalization of Sylow's third theorem,
J. Algebra, 115 (1988), 414-430.

[Frob95]
F. G. Frobenius,
Verallgemeinerung des Sylow'schen Satzes,
Sitzungsberichte der K\"onigl. Preu\SS. Akad. der Wissenschaften (Berlin)
(1895), 981-993.

[Frob03]
F. G. Frobenius, \"Uber einen Fundamentalsatz der Gruppentheorie,
Sitzungsberichte der K\"onigl. Preu\SS. Akad. der Wissenschaften (Berlin)
(1903), 987-991.

[GRV12]
C. Gordon, F. Rodriguez-Villegas,
On the divisibility of $\#\Hom(\Gamma, G)$ by $|G|$,
J. Algebra, 350:1 (2012), \hbox{300-307}.
\arXiv 1105.6066

[Hall36]
Ph. Hall,
On a theorem of Frobenius,
Proc. London Math. Soc. 40 (1936), 468-501.



[Isaa70]
I. M. Isaacs,
Systems of equations and generalized characters in groups,
Canad. J. Math., 22 (1970), \hbox{1040-1046}.

[Iwa82]
S. Iwasaki,
A note on the $n$th roots ratio of a subgroup of a finite group,
J. Algebra, 78:2 (1982), 460-474.

[KM14]
A. A. Klyachko, A. A. Mkrtchyan,
How many tuples of group elements have a given property?
With an appendix by Dmitrii V. Trushin,
Intern. J. of Algebra and Comp. 24:4 (2014), 413-428.
\arXiv 1205.2824

[KM17]
A. A. Klyachko, A. A. Mkrtchyan,
Strange divisibility in groups and rings,
Arch. Math. 108:5 (2017), 441-451.
\arXiv 1506.08967


[Kula38]
A. Kulakoff,
Einige Bemerkungen zur Arbeit: ``On a theorem of Frobenius" von P. Hall,
Mat. Sb.,  3(45):2 (1938), 403-405.

[SaAs07]
J. Sato, T. Asai,
On the $n$-th roots of a double coset of a finite group,
J. School Sci. Eng., Kinki Univ., 43 (2007), 1-4.

[Sehg62]
S. K. Sehgal,
On P. Hall's generalisation of a theorem of Frobenius,
Proc. Glasgow Math. Assoc.,  5 (1962), 97-100.

[Solo69]
L. Solomon,
The solution of equations in groups,
Arch. Math., 20:3 (1969), 241-247.


[Yosh93]
T. Yoshida,
$|\Hom(A, G)|$,
Journal of Algebra, 156:1 (1993), 125-156.

\endREFERENCES

\end